\documentclass[leqno,11pt]{amsart}

\usepackage{amssymb}
\usepackage{amsmath}
\usepackage{enumerate}
\usepackage{amsfonts}
\usepackage{color}
\usepackage[colorlinks,citecolor=blue,urlcolor=blue]{hyperref}	
\usepackage[final]{graphicx}

\newtheorem{teo}[equation]{Theorem}

\newtheorem{lema}[equation]{Lemma}
\newtheorem{propo}[equation]{Proposition}
\newtheorem{remark}[equation]{Remark}

\numberwithin{equation}{section}

\newcommand{\al}{\alpha}
\renewcommand{\a}{\alpha}
\renewcommand{\sl}{\slash}
\newcommand{\SSS}{\Omega}
\newcommand{\SSSS}{\Omega'}

\newcommand{\xt}{\frac{x}{2t}}
\renewcommand{\d}{\partial}
\newcommand{\dx}{\frac{\d}{\d x}}

\newcommand{\dtt}{\frac{dt}{\sqrt{t}}}
\newcommand{\dttt}{\frac{dt}{t}}
\renewcommand{\xt}{\frac{x}{2t}}

\newcommand{\C}{\mathbb{C}}

\newcommand{\mc}{\mathcal}

\newcommand{\supp}{\mathrm{supp}}
\newcommand{\N}{\mathbb{N}}

\newcommand{\e}{\varepsilon}

\newcommand{\8}{\infty}

\newcommand{\LL}{L^1(X)}

\newcommand{\Hardd}{\wt{H}^1(X)}
\newcommand{\Harddat}{\wt{H}^1_{at}(X)}

\renewcommand{\(}{\left(}
\renewcommand{\)}{\right)}

\newcommand{\wt}{\widetilde}

\newcommand{\ii}{\int_0^{\infty}}
\newcommand{\iie}{\int_\e^{\e^{-1}}}
\newcommand{\pssii}{\omega}

\begin{document}

\title[Characterization of Hardy spaces for the Laguerre expansion]
{Riesz transform characterization of $H^1$ spaces associated with
certain  Laguerre expansions}

\begin{abstract}
For $\alpha>0$ we consider the system $\psi_k^{(\alpha-1)\slash 2}(x)$  of
the Laguerre functions which are eigenfunctions of the
differential operator $Lf = -\frac{d^2}{dx^2}f -
\frac{\alpha}{x}\frac{d}{dx}f + x^2 f$. We define an atomic Hardy
space $H^1_{at}(X)$, which is a subspace of $L^1((0,\8), x^\alpha \,
dx)$. Then we prove that the space $H^1_{at}(X)$ is also characterized by the Riesz
transform $R f =\sqrt{\pi} \dx L^{-1/2} f$ in the sense that $f\in
H^1_{at}(X)$ if and only if $f, Rf \in L^1((0,\8),x^\alpha \, dx)$.
\end{abstract}

\author[Marcin Preisner ]{ Marcin Preisner}
\address{Instytut Matematyczny\\
Uniwersytet Wroc\l awski\\
50-384 Wroc\l aw, Pl. Grunwaldzki 2/4\\
Poland} \email{preisner@math.uni.wroc.pl}

\thanks{
The paper is a part of the author's Ph.D. thesis. The research was partially supported by Polish funds for sciences, grants: N N201 397137 and N N201 412639, MNiSW}

\maketitle

\begin{flushleft}
\end{flushleft}

\section{Background and main result}\label{sect1}
For fixed $\a >0$ let $X$ denote the space $(0,\8)$ with the
measure $d\mu(x) = x^\a dx$. The space $X$ equipped with the
Euclidean distance $d(x,y)=|x-y|$ is a space of homogeneous type
in the sense of Coifman-Weiss \cite{CW2}. On $L^2(X)$ we consider
the orthogonal system of the Laguerre functions
$\{\psi_k^{(\a-1)\slash 2}(x)\}_{k=0}^\infty$,
$$\psi_k^{(\a-1)\slash 2}(x) = \(\frac{2k!}{\Gamma(k+\a \sl 2+1\sl 2)}\)^{1\sl 2} L_k^{(\a-1)\slash 2}
(x^2)e^{-x^2\sl 2},$$
where $L_k^\a$ is the $k$-th Laguerre
polynomial. Each $\psi_k^{(\a-1)\slash 2}$ is an eigenfunction of
the Laguerre operator
\begin{equation*}
L f (x)= -\frac{d^2}{dx^2}f(x) - \frac{\a}{x}\frac{d}{dx}f(x) + x^2
f(x),
\end{equation*}
where the corresponding eigenvalue is $\beta_k=4k+\a+1$. Let
 $$ T_tf=\sum_{k=0}^\infty \exp(-t\beta_k) \langle f, \psi_k^{(\a-1)\slash 2}\rangle \psi_k^{(\a-1)\slash 2}$$
 be the semigroup of the self-adjoint linear operator on $L^2(X)$ generated by $-L$,
where  $\mathcal D(L)=\{ f\in L^2(X): \sum_{k} \beta_k^2 |\langle
f,\psi_k^{(\a-1)\slash 2}\rangle |^2 <\infty \}$ is the domain of
$L$.

It is well known (see e.g. \cite{Leb}, \cite{NS}) that $T_t$ has the integral
representation, i.e.,
\begin{equation}\label{semi0}
T_t f (x) = \ii T_t(x,y) f(y) d\mu (y),
\end{equation}
where
\begin{equation}\label{semi}
T_t(x,y) = \frac{2e^{-2t}(xy)^{-(\a-1)\sl 2}}{1-e^{-4t}}
\exp\(-\frac{1}{2}\frac{1+e^{-4t}}{1-e^{-4t}}(x^2+y^2)\) I_{(\a-1)\sl 2} \(\frac{2e^{-2t}}{1-e^{-4t}}xy\).
\end{equation}
Here $I_\nu$ denotes the Bessel function of the second kind.
The operators (\ref{semi0}) define strongly
continuous semigroups of contractions on every $L^p(X)$, $1\leq
p<\infty$.

Through this paper we shall use the following notation: for an interval $I \subseteq (0,\8)$ we will denote by $|I|$ its Euclidean diameter, $B(x,r) =\{y \in X \, : \, |x-y|<r\}$, and $\chi_A$ will be the characteristic function of the set $A$. We define the auxiliary  function
\begin{equation}\label{rho}
\rho(y)=\chi_{(0,1)}(y)+\frac{1}{y}\chi_{[1,\8)}(y).
\end{equation}

{\bf Definition.} A function $a$ is called an $H^1(X)-${\it atom} if there exists an interval $I = B(y_0,r) \subseteq (0,\8)$ such that:
\begin{enumerate}
\item [(i)]$\supp (a) \subseteq I \ \text{ and }\  r\leq \rho(y_0)$,
\item [(ii)]$\|a\|_\8 \leq \mu(I)^{-1}$,
\item [(iii)]if $r \leq \rho(y_0)/4$, then $\ii a(x) d\mu (x) = 0$.
\end{enumerate}
We say that an $L^1(X)$-function $f$ belongs to $H^1_{at}(X)$ if and
only if there exist sequences $\{a_j, \lambda_j\}_{j=1}^\8$ such that $f=\sum_{j=1}^\8 \lambda_j a_j$, each $a_j$ is an
$H^1(X)$-atom, $\lambda_j\in \C$, and $\sum_{j=1}^\8 |\lambda_j|
<\8$. The space $H_{at}^1(X)$ is a Banach space with the norm

\begin{equation*}\label{HarLagNorm}
\|f\|_{H^1_{at}(X)} = \inf \sum_{j=1}^{\8} |\lambda_j|,
\end{equation*}
where the infimum is taken over all representations $f=\sum_{j=1}^{\8} \lambda_j a_j$.

Let $\delta=\frac{d}{dx}+x$,
$\delta^*=-\frac{d}{dx}+x-\frac{\a}{x}$. Then
$L=(\a+1)I+\delta^*\delta$, $\delta \psi_k^{(\a-1)\slash 2}=
-2\sqrt{k}x\psi_{k-1}^{(\a+1)\slash 2}$.  The Riesz transform
$R_L$, originally defined on $L^2(X)$ (see, e.g., \cite{NS2}, \cite{NS}) by the formula
$$ R_L f=\sqrt{\pi}\delta L^{-1\slash 2} f=-\sum_{k=1}^\infty
\Big(\frac{4k\pi}{4k+\a+1}\Big)^{1\slash 2} \langle
f,\psi_k^{(\a-1)\slash 2}\rangle x \psi_{k-1}^{(\a+1)\slash 2},$$
turns out to be the principal value
singular integral operator
 $$ R_L f(x)=\lim_{\varepsilon \to 0}
 \int_{0, \, |x-y|>\varepsilon }^\infty  R_L (x,y)f(y) d\mu(y),$$
 with the kernel
 $$ R_L (x,y)=\int_0^\infty \Big(\frac{\partial}{\partial x}+x\Big)T_t(x,y)\frac{dt}{\sqrt{t}}.$$
 Since the kernel
 $$\Gamma (x,y)=\int_0^\infty xT_t(x,y)\frac{dt}{\sqrt{ t}}$$
 satisfies $\sup_{y>0}\int|\Gamma (x,y)|d\mu(x)<\infty$, it defines a bounded linear operator on
 $L^1(X)$.
 Hence, for our purposes, we restrict our consideration to the Riesz transform $Rf=\sqrt{\pi}\frac{d}{dx}L^{-1\slash 2}f$.
Clearly, $R$ is a principal value
 singular integral operator with the kernel
\begin{equation}\label{Rieszaaa}
R(x,y)=\ii \dx T_t(x,y)\dtt.
\end{equation}
The action of $R$ on $L^1(X)$-functions is well-defined in the sense of distributions (see Section \ref{rieszkerlag} for details).

The main goal of this paper is to prove the following theorem.
\begin{teo} \label{MainThm}
A function $f \in \LL$ belongs to the Hardy space $H^1_{at}(X)$ if and only if $Rf$ belongs to $\LL$. Moreover, the corresponding norms are equivalent, i.e.,
\begin{equation}\label{mainineq}
C^{-1} \|f\|_{H^1_{at}(X)} \leq \|f\|_{\LL}+\|Rf\|_{\LL} \leq C
\|f\|_{H^1_{at}(X)}.
\end{equation}
\end{teo}

The main idea of the proof is to compare the kernel $R(x,y)$ with
 kernels of appropriately scaled   local Riesz transforms related to the
 Bessel operator $\widetilde L f(x)=-f''(x) -\frac{\a}{x}f'(x)$, where
 the scale of localization is adapted to the auxiliary function
 $\rho(y)$.
 To do this we consider the Bessel semigroup:
\begin{align}\label{kerbess}
&\wt{T}_t f (x)= \ii \wt{T}_t (x,y) f(y) d\mu(y), \cr
&\wt{T}_t(x,y) = (2t)^{-1} \exp\left(-\frac{x^2+y^2}{4t}\right) I_{(\a-1)/2} \left(\frac{xy}{2t}\right) (xy)^{-(\a-1)/2}
\end{align}
and observe that for small $t$ the kernel \eqref{kerbess} is close to the kernel
\eqref{semi}. Thanks to this, $R(x,y)$ is
comparable to $\wt R (x,y)$  after some suitable localization
defined by the function $\rho$, where $\wt R(x,y)$ denotes the
Riesz transform kernel in the Bessel setting. This
requires  a precise computation of constants appearing in singular
parts of the kernels (see Propositions \ref{lem2.2} and \ref{lem3.3}).
The next step is to use results of Betancor, Dziuba\'nski, and Torrea \cite{BDT}, who give characterizations of a "global"
Hardy space for the Bessel operator, to define and describe
local Hardy spaces for $\widetilde L$. Having all these
prepared we prove the theorem.

We would like to remark that the Hardy space $H^1_{at}(X)$ we consider here is also characterized by means of the maximal function:
$$M f(x) = \sup_{t>0} |T_t f(x)|,$$
that is, the norm $\|f\|_{H^1_{at}(X)}$ is comparable with $\|Mf\|_{\LL}$. For details we refer the reader to \cite{Dz3}.

There are other expansions based on the Laguerre functions for which Hardy spaces were investigated. For example, when $\a >-1$, systems $\{\varphi_n^\a\}_{n=0}^\8$ and $\{\mathcal{L}_n^\a\}_{n=0}^\8$, where
$$\varphi_n^\a (x) = c_{n,\a} e^{-x^2/2}x^{\a+1/2}L_n^\a(x^2), \qquad \mathcal{L}_n^\a (x) = c_{n,\a}e^{-x/2}x^{\a/2} L_n^\a(x),$$
are orthogonal on $L^2((0,\8), dx)$. These systems are related to operators
$$\mathbf L _\a = -\frac{d^2}{dx^2} + x^2 +\frac{1}{x^2}(\a^2-\frac{1}{4}), \qquad \overline{L}_\a= -x\frac{d^2}{dx^2} -\frac{d}{dx}+\frac{x}{4}+\frac{\a^2}{4x},$$
respectively. In \cite{DGT} and \cite{Dz2} the authors proved that the Hardy spaces associated with $\{\varphi_n^\a\}_{n=0}^\8$ and $\{\mathcal{L}_n^\a\}_{n=0}^\8$ are characterized by: the maximal functions, the Riesz transforms, and certain atomic decompositions. Moreover, in \cite{Dz1} the author obtained an atomic description of the Hardy space originally defined by the maximal function related to the system
$$\ell _n^\a (x) = c_{n,\a} L_n^\a(x) e^{-x/2}, \quad n=0,1,..., \  \ \text{on } \ L^2((0,\8), x^\a dx).$$
The functions $\ell_n^\a$ are eigenfunctions of the operator
$$\mathbb L_\a = -x\frac{d^2}{dx^2}-(\a+1)\frac{d}{dx}+\frac{x}{4}.$$

Finally, we would like to note that the system $\{\psi_n^{(\a-1)/2}\}_{n=0}^\8$ we consider in the present paper is well-defined and orthogonal on $L^2(X)$ for $\a >-1$. However, the case $-1<\a\leq 0$ is not included in our investigations.

The paper is organized as follows. In Section \ref{HarBess} we present a singular integral characterization of local Hardy
spaces associated with the Bessel operator $\widetilde L$. Section
 \ref{rieszkerlag} is devoted to stating detailed estimates for $R(x,y)$ and
 proving some auxiliary results. The proof of the main theorem is given in Section
 \ref{Mainproof}.  In Section \ref{aux} we present proofs of estimates of the kernels $\wt R(x,y)$ and $R(x,y)$ stated in Propositions \ref{lem2.2} and \ref{lem3.3}.

\section{Hardy spaces for the Bessel operator} \label{HarBess}
\subsection{Global Hardy space}\label{globalBess}
Hardy spaces $\wt{H}^1(X)$ related to Bessel operators were studied in \cite{BDT}.

{\bf Definition.} We call a function $a$ an $\Hardd-${\it atom} if there is an interval $I \subset (0,\8)$ such that:
\begin{enumerate}
\item [(i)]$\supp (a) \subseteq I$,
\item [(ii)]$\|a\|_\8 \leq \mu(I)^{-1}$,
\item [(iii)]$\ii a(x) d\mu (x) = 0$.
\end{enumerate}
We define the space $\Harddat$ in the same way as $H^1_{at}(X)$ in Section \ref{sect1}.

The singular integral kernel of the Riesz transform $\wt{R}$ is defined by
\begin{equation*}\label{RieszBess}
\wt{R}(x,y) =\ii \dx \wt{T}_t (x,y) \dtt, \qquad \text{where } \ x\neq y.
\end{equation*}
Before giving a distributional sense of $\wt R f$ for $f \in \LL$ we recall results from \cite{BDT}.
\begin{teo} \label{MainBess}
For $f\in \LL$ the following conditions are equivalent:
\begin{enumerate}
\item [(i)] $f \in \wt H^1_{at}(X),$
\item [(ii)] $\wt R f \in \LL$,
\item [(iii)] $\sup_{t>0}|\wt T_t f| \in \LL.$
\end{enumerate}
Moreover,
$$\|f\|_{\Harddat} \sim \(\|f\|_{\LL} + \|\wt{R} f\|_{\LL}\) \sim \Big\|\sup_{t>0}|\wt T_tf|\Big\|_{\LL}.$$
\end{teo}

For a function $f$ defined on $(0,\8)$ and $y>0$ we denote $f_y(x) = y^{-\a-1} f(x/y)$. Let
\begin{equation}\label{const}
A= A(\a)= -\frac{2 \Gamma (1+\a/2)}{ \Gamma ((1+\a)/2)}=-\frac{2\gamma_1}{\gamma_2}, \qquad B=B(\a)= - \frac{\a+1}{\sqrt{\pi}}.
\end{equation}
The following proposition will play a crucial role in our investigations.
\begin{propo}\label{lem2.2}
Let $A,B$ be as in \eqref{const}. Then for $x\neq y$ we have
$$
\wt{R} (x,y) = \frac{A-B}{x^{\a+1} + y^{\a+1}} + \frac{B}{x^{\a+1} - y^{\a+1}} + h_y(x),
$$
where
\begin{equation}\label{prop_h}
h \in \LL \quad  \text{ and  } \quad |h(x) + A-2B| \leq Cx \quad \text{  for  }x\leq 1/2.
\end{equation}
\end{propo}

The proof of Proposition \ref{lem2.2} is postponed until Section \ref{proof2.2}. To give a precise definition of $\wt R$ on $\LL$ we need a suitable space of test functions. One of possible choices is
$$\SSS (X) = \left\{ \pssii \in C^{1} (0,\8) \ \big| \ \ \|\pssii\|_\8, \ \Big\|\frac{\pssii(x)}{x}\Big\|_{L^1(X)}, \ \|x\pssii'(x)\|_\8<\8 \right\}$$
with the topology defined by the semi-norms $\gamma_i$, $i=1,2,3$, where,
$$\gamma_1(\pssii)= \|\pssii\|_\8, \quad \gamma_2(\pssii)=\Big\|\frac{\pssii(x)}{x}\Big\|_{\LL}, \quad \gamma_3(\pssii)=\|x\pssii'(x)\|_\8.$$
Denote by $\SSSS(X)$ the dual space.

The space $f \in \LL$ is contained in $\SSSS(X)$ in the natural sense, i.e., if $f\in \LL$, then
$$<f,\pssii > = \ii f \pssii \, d\mu, \qquad \pssii \in \SSS(X).$$

Next, for $f\in \LL,\  \pssii \in \SSS (X)$, we define
\begin{equation}\label{distr2}
\langle \wt{R}f, \pssii \rangle = \langle f, \wt{R}^* \pssii \rangle, \qquad \wt{R}^* \pssii (y) = \lim_{\e \to 0} \int_{|x-y|>\e}  \wt{R} (x,y)\pssii (x) d\mu(x).
\end{equation}
Alternatively, we define the Riesz transform as follows:
\begin{equation}\label{distr}
\langle \wt{\mc R}f, \pssii \rangle =  \langle f, \wt{\mc R}^* \pssii \rangle, \qquad  \wt{\mc R}^* \pssii (y)=  \lim_{\e \to 0} \iie \int_X \dx \wt T_t(x,y) \pssii (x) d\mu(x) \dtt.
\end{equation}

\begin{propo}\label{lllppp}
For $\pssii \in \SSS (X)$ and $y>0$ we have
$\wt R^* \pssii (y)= \wt {\mc R}^* \pssii(y)$.
Moreover,
\begin{equation*}
\| \wt R^* \pssii \|_\8 \leq C\(\|\pssii(x)\|_\8+\|x\pssii'(x)\|_\8+\Big\|\frac{\pssii(x)}{x}\Big\|_{L^1(X)}\) .
\end{equation*}
\end{propo}
The proof can be deduced from \eqref{kerbess} and Proposition \ref{lem2.2}. We will not go into details here. However, we would like to notice that from Proposition \ref{lllppp} it is easily seen, that the Riesz transform $\wt{\mc R}$ is the same as the one defined by the spectral theorem (see, e.g., \cite{NS2}--\cite{NS}).

\subsection{Local Hardy spaces}\label{localBess}
Fix a non-negative function $\phi \in C_c^\8 (-2,2)$ such that $\phi(x) = 1$ for $|x|\leq 3/2$. Similarly to the classical case, for $m>0$ we define {\it scaled  local Riesz transforms} $\wt r^{m}$ for $f\in \LL$, $\pssii \in \SSS(X)$  as follows:
$$\langle \wt{r}^m f,\pssii \rangle = \langle f,(\wt{r}^{m})^* \pssii\rangle , \quad  (\wt{r}^{m})^* \pssii (y) = \lim_{\e\to 0} \int^\8_{0,|x-y|>\e} \wt{R} (x,y) \phi\(\frac{x-y}{m}\) \pssii(x) d\mu (x). $$
As in the global case these operators are well-defined and
\begin{equation}\label{cont}
\|(\wt r^{m})^* \pssii\|_\8 <\8.
\end{equation}

For an interval $I=B(y,r) \subseteq X$ and $k>0$ let $kI = B(y,kr) \subseteq X$.

\begin{lema}\label{L2bound}
The operators $\wt{r}^m$ are bounded on $L^2(X)$ with norm-operator bounds independent of $m$.
\end{lema}
\begin{proof} Because of the dilatation structure (see \eqref{dill}) it is enough to prove the lemma in the case $m=1$. Assume additionally that $\supp f \subseteq I = B(y_0, 1)$. Then $\wt r ^1 f(x) = 0$ for $x \notin 3I$. Also
$$\|\wt r ^1 f\|_{L^2(X \cap 3I)} \leq \|(\wt r ^1- \wt R)f\|_{L^2(X \cap 3I)} + \|\wt R f\|_{L^2(X)}.$$
It is well known that $\|\wt R f\|_{L^2(X)} \leq C\|f\|_{L^2(X)}$ (see \cite{MS}). Moreover,
\begin{align}\label{est15}
|\wt R(x,y)|\chi_{\{|x-y|>3/2\}} \leq C(xy)^{-\a/2} + |h_y(x)| \leq C (xy)^{-\a/2} + |h_1(x,y)|+|h_2(x,y)|,
\end{align}
where $$h_1(x,y) = y^{-\a-1}  (h-D\chi_{(0,1)})(x/y), \quad  		h_2(x,y)=D(\chi _{(0,1)})_y(x) = Dy^{-\a-1}  \chi_{(x,\8)}(y).$$
Here $D=A-2B$. We claim that
$$\|(\wt r^1 - \wt R)f\|_{L^2(X\cap 3I)} \leq C \|f\|_{L^2(X)}.$$
To prove this we consider the three summands from \eqref{est15} separately. By the Cauchy-Schwarz inequality we get
\begin{align*}
\Big\|\int_I (xy)^{-\a/2} f(y) d\mu (y)\Big\|^2_{L^2(X\cap 3I)} \leq C \int_{3I} \|f\|^2_{L^2(X)} dx \leq C\|f\|^2_{L^2(X)}.
\end{align*}
From \eqref{prop_h} we deduce
$$\sup_{y>0}\ii |h_1(x,y)|d\mu (x) + \sup_{x>0}\ii |h_1(x,y)|d\mu(y) <\8.$$
Thus the operator with the kernel $h_1(x,y)$ is bounded on every $L^p(X)$, $1\leq p \leq \8$. The part which contains $h_2$ is bounded on $L^2(X)$ due to the Hardy inequality (see, e.g. \cite{Bennett}, p. 124).

To omit the assumption $\supp f \subseteq B(y_0, 1)$ let us notice
$$\|\wt r^1 f\|_{L^2(X)} \leq \sum_{j=1}^\8 \|\wt r^1 (f \cdot \chi_{(j-1,j)})\|_{L^2(X)} \leq C \sum_{j=1}^\8 \|f \cdot \chi_{(j-1,j)}\|_{L^2(X)} = C \|f\|_{L^2(X)}.$$
\end{proof}

The local Hardy space $\wt{h}^{1,m} (X)$ is a subspace of $\LL$ consisting of functions $f$ for which $\wt{r}^m f \in \LL$. In order to state atomic characterization of $\wt{h}^{1,m}(X)$ we call a function $a$ an $\wt{h}^{1,m}(X)$-atom, when there exists an interval $I=B(y_0,r) \subset (0,\8)$ such that
\begin{enumerate}
\item[(i)] $\supp (a) \subseteq I\ \text{ and } \ r\leq m$,
\item[(ii)] $\|a\|_\8 \leq \mu(I)^{-1}$,
\item[(iii)] if $r \leq m/4$, then $\ii a(x) d\mu (x) = 0$.
\end{enumerate}

\begin{teo}\label{localdecomp}
Assume that $f\in \LL$. Then $\wt{r}^m f \in \LL$  if and only if there exist sequences $\lambda_j \in \C$ and $\wt h^{1,m}(X)$-atoms $a_j$, such that $f=\sum_{j=1}^\8 \lambda_j a_j$, where $\sum_{j=1}^\8 |\lambda_j| < \8$. Moreover, we can choose $\{\lambda_j\}_{j}, \{a_j\}_{j}$, such that
$$C^{-1}\sum_{j=1}^\8 |\lambda_j| \leq  \|f\|_{\LL}+ \|\wt{r}^{m}f\|_{\LL} \leq C \sum_{j=1}^\8 |\lambda_j|,$$
where $C	$ is independent of $m>0$.
\end{teo}
\begin{remark}\label{remark}
Assume in addition $\supp (f) \subseteq I=B(y_0, m)$. Then, in the above decomposition, one can take atoms with supports contained in $3 I$.
\end{remark}
\begin{proof}
The proof is similar to the classical case. For the reader's convenience we provide some details. Without loss of generality we may assume that $m=1$. The operator $\wt r ^{1}$ is continuous from $\LL$ to $\SSS '(X)$ (see \eqref{cont}), so the first implication will be proved when we have obtained
\begin{equation}\label{fgj}
\|\wt r^{1} a\|_{\LL} \leq C
\end{equation}
for every $\wt h^{1,1}(X)$-atom $a$. Notice, that the weak-type $(1,1)$ bounds of $\wt r^1$ also reduces the proof to \eqref{fgj}, as it was pointed out to us by the referee. Assume then, that $a$ is an $\wt h ^{1,1}(X)$-atom supported by an interval $I=B(y_0,r)$. Note that $\wt r^1 a(x) =0$ on $(9I)^c$. Consider first the case where $r>1/4$. Recall that $\mu$ has the doubling property. By the Cauchy-Schwarz inequality and Lemma \ref{L2bound} we get
$$\|\wt r^1 a\|_{L^1(X \cap 9I)} \leq \mu (9I)^{1/2} \|\wt r^1 a\|_{L^2(X)}  \leq C \mu (I)^{1/2} \|a\|_{L^2(X)}  \leq  C.$$
If $r<1/4$ then $a$ is an $\wt H ^1(X)-$atom, so by Theorem \ref{MainBess} it follows that $\|\wt Ra\|_{\LL} \leq C$. Therefore $\| \wt r^1 a\|_{\LL} \leq C+ \|(\wt R-\wt r^1) a\|_{\LL}$. Because of the cancelation condition we have
\begin{equation*}
(\wt R - \wt r^1) a(x) = \int \( \wt R (x,y)(1- \phi(x-y)) - \wt R (x,y_0) (1-\phi(x-y_0)) \) a(y) d\mu (y).
\end{equation*}
Thus it is enough to verify the estimate
\begin{equation}\label{estXI}
\sup _{y\in I} \ii \big| \wt R (x,y) (1-\phi(x-y)) - \wt R (x,y_0) (1-\phi(x-y_0)) \big| d\mu(x) = \sup _{y\in I} \ii \Xi(x,y) d\mu(x) \leq C.
\end{equation}
Fix $y \in I$. From Proposition \ref{lem2.2} one obtains
\begin{align}\nonumber
&\Xi(x,y) = 0 \ &&\text{for} \ |x-y_0|\in (0,1)\\ \label{TUV}
&\Xi(x,y) \leq C \, x^{-\a} + |h_y(x)|+ |h_{y_0}(x)| \ &&\text{for} \ |x-y_0|\in (1, 3)\\
&\Xi(x,y) \leq C \, |x-y_0|^{-2} x^{-\a} + |h_y(x)|+ |h_{y_0}(x)|  &&\text{for} \ |x-y_0|\in (3,\8), \nonumber
\end{align}
where in the last inequality we have used that $\phi(x-y)=\phi(x-y_0)=0$ and the mean-value theorem.
From \eqref{TUV} we get \eqref{estXI} and $\|(\wt R-\wt r^1) a\|_{\LL} \leq C$. This ends the proof of \eqref{fgj}.

For the converse, assume that $f, \wt r ^{1} f \in \LL$ and, in addition, $\supp f \subseteq I = B(y_0, 1).$ Fix $\xi = \mu(I)^{-1} \int_I f d\mu$, $g=f-\xi \chi_I$. We have
\begin{equation}\label{dec}
\|\wt R g \|_{\LL} \leq \|\wt r^{1} f\|_{\LL} + \|\xi \wt r^{1}(\chi_I)\|_{\LL} + \|(\wt R - \wt r^{1})g\|_{\LL}.
\end{equation}
By using the first part of the proof we deduce that $\|\xi \wt r^1 \chi_I\|_{\LL} \leq C\|f\|_{\LL}$. Note that $\supp \,  g \subseteq I$, $ \int g \, d\mu= 0$, so \eqref{estXI} implies $\|(\wt R - \wt r^1) g\|_{\LL}\leq C\|g\|_{\LL} \leq C\|f\|_{\LL}$. Therefore, by Theorem \ref{MainBess}, there exist $\wt H^1(X)$-atoms $a_j$ ($j=1,...$) such that
$$f-\xi \chi_I = g= \sum_{j=1}^\8 \lambda_j a_j.$$
Moreover $\sum_{j=1}^\8 |\lambda_j| \leq \|f\|_{\LL} +\|\wt r^{1}f\|_{\LL}$. Denote $\lambda_0 = \int_I f d\mu, \ \ a_0 = \mu(I)^{-1} \chi_I$ and fix $\psi_I \in C_c^\8\(\frac{4}{3}I\)$ satisfying $\psi_I \equiv 1$ on $I$ and $\|\psi_I\|_{\8} \leq C$. We have obtained
\begin{equation}\label{dec1}
f= f \psi_I=\sum_{j=0}^\8 \lambda_j (\psi_I a_j).
\end{equation}
It remains to show that each $\psi_I a_j$ can be written in the following form: $\psi_I a_j = \sum_{i=1}^{N_j} \kappa_{i,j} b_{i,j}$, where $b_{i,j}$ are $\wt h ^{1,1}(X)$-atoms supported in $3I$ and $\sum_{i=1}^{N_j} |\kappa_{i,j}| \leq C$, where $C>0$ is independent of $j$. For $j=0$ the claim is clear. Fix $j\geq 1$ and suppose that $\supp \, a_j \subseteq J = B(z_0,r)$. Obviously, if~$(\frac{4}{3} I) \cap J = \emptyset$ then $\psi_I a_j =0$. Moreover, if $r>1/4$ then $\psi_I a_j = \kappa b$, where $b$ is an $\wt h ^{1,1}(X)$-atom and $|\kappa| \leq C$. So, suppose that $(\frac{4}{3} I) \cap J \neq \emptyset$ and $r<1/4$. Under these assumptions we write
\begin{align*}
\psi_I(x) a_j(x) = &\(\psi_I(x) a_j(x)-\sigma \mu(2J)^{-1} \chi_{2J}(x)\) \cr
&+ \sum_{i=1}^{N-1} \sigma \((\mu(2^i J))^{-1} \chi_{2^i J} (x) - (\mu(2^{i+1} J))^{-1} \chi_{2^{i+1} J} (x)\)\cr
&+\sigma (\mu(2^{N} J))^{-1} \chi_{2^{N} J} (x),
\end{align*}
where $\sigma = \ii a_j(z)(\psi_I(z)-\psi_I(z_0))d\mu(z)$ and $N$ is such that $2^{-N-1}\leq r < 2^{-N}$. One can check that this is the required decomposition, since $|\sigma|\leq C r$. Let us note that we have just proved Remark \ref{remark}.

To deal with the general case we take a smooth partition of unity $\{ \psi_j \}_{j=1}^\8 \subseteq C^\8 (0,\8)$, i.e.
$$\sum_{j=1}^\8 \psi_j(x) = \chi_{(0,\8)}(x), \quad 0\leq \psi_j \leq 1, \quad \supp \, \psi_j\subseteq I_j = B\(y_j, 1\), \quad \sup_{j\in \N} \|\psi '_j\|_\8 \leq C. $$
Consider
$$g_j = \wt r^1(\psi_j f)- \psi_j \wt r^1(f). $$
Obviously, $\supp \, g_j \subseteq 3I_j$ and for $x \in 3I_j$ we have
\begin{align}\label{ttt}
|g_j(x)| &= \Big|\ii \wt R(x,y) \phi(x-y) f(y)(\psi_j(y)-\psi_j(x)) d\mu(y)\Big|\cr
&\leq C \ii |\wt R(x,y)| \chi_{\{|x-y|\leq 2\}}|f(y)| |x-y| d\mu(y).
\end{align}
Moreover, from Proposition \ref{lem2.2} we have
\begin{equation}\label{sss}
\sup_{y>0}\int_{|x-y|\leq 2}|\wt R (x,y)||x-y|d\mu(x) \leq C.
\end{equation}
From \eqref{ttt} and \eqref{sss} we deduce
$$\|g_j\|_{L^1(X)} \leq \|f\|_{L^1(X\cap 5I_j)}.$$
Therefore
\begin{equation}\label{ppp}
\sum_{j=1}^\8 \|\wt r^{1} (\psi_j f) \|_{\LL} \leq \sum_{j=1}^\8 \( \|\psi_j \wt r^{1} ( f) \|_{\LL} +  \| g_j \|_{\LL}\)
\leq C\(\|f\|_{\LL} + \|\wt r^1(f)\|_{\LL}\).
\end{equation}

By using \eqref{dec1} and the subsequent remark for each $\psi_j f$ we get the decomposition $\psi_j f = \sum_{k} \lambda_k^j a_k^j$, where $a_k^j$ are $\wt h^{1,1}(X)$-atoms and
\begin{equation}\label{ggg}
\sum_k |\lambda_k^j| \leq C\(\|\psi_j f\|_{\LL} + \|\wt r^1 (\psi_j f)\|_{\LL}\).
\end{equation}
The proof is completed by noticing that
$$f = \sum_{j,k} \lambda_k^j a_k^j,$$
where
$$ \sum_{j,k} |\lambda_k^j| \leq C\( \|f\|_{\LL} + \|\wt r^1 f\|_{\LL}\)$$
is guaranteed by \eqref{ppp} and \eqref{ggg}.
\end{proof}

\section{The Riesz kernel for the Laguerre expansion}\label{rieszkerlag}

Let $\phi$ be the function defined in Section \ref{localBess} and $\rho$ be as in \eqref{rho}. The following proposition gives an essential information about the kernel of the Riesz transform associated with the Laguerre expansion.

\begin{propo}\label{lem3.3} Let $A$ and $B$ be as in \eqref{const}. The kernel $R(x,y)$ can be written in the form
$$R(x,y) =\phi\(\frac{x-y}{\rho(y)}\) \( \frac{B}{x^{\a+1}-y^{\a+1}} + \frac{A-B}{x^{\a+1}+y^{\a+1}} \) + g(x,y),$$
where
\begin{equation}\label{integr}
\sup_{y>0} \ii |g(x,y)| d\mu(x) < \8.
\end{equation}
\end{propo}

The proof of Proposition \ref{lem3.3} is a quite lengthy analysis. We provide details in Section \ref{proof3.3}.

For $f\in \LL,\  \pssii \in \SSS (X)$, we define the Riesz transform $Rf$ as follows
\begin{align*}
&\langle Rf, \pssii \rangle = \langle f, R^* \pssii \rangle, \qquad R^* \pssii (y) = \lim_{\e \to 0} \int_{|x-y|>\e}  R (x,y)\pssii (x) d\mu(x),\cr
\end{align*}

One can easily check using Proposition \ref{lem3.3} that this limit exists and
\begin{align} \label{cont2}
&\| R^* \pssii \|_\8 \leq C\(\|\pssii(x)\|_\8+\|x\pssii'(x)\|_\8+\Big\|\frac{\pssii(x)}{x}\Big\|_{L^1(X)}\) .
\end{align}

Denote by $G$ the operator with the kernel $g(x,y)$. Obviously, by \eqref{integr}, $G$ is bounded on $\LL$. In the proof of Theorem \ref{MainThm} we will need the following lemma.
\begin{lema}\label{comut}
Let $z\in (0,\8)$, $f \in \LL$, $I = B(z, \rho(z))$, and $\eta \in C^\8\(0,\8\)$ satisfies $0 \leq \eta \leq 1$, $\supp \ \eta \subset I$, $\|\eta '\|_\8 \leq C_1 \rho(z)^{-1}$. Then
$$\|R(\eta f) -\eta ((R-G)f)\|_{\LL} \leq C \|f\|_{L^1(X \cap 4I)},$$
with a constant $C$ which depends on $C_1$, but it is independent of $z\in (0,\8)$ and $f \in \LL$.
\end{lema}
\begin{proof} Note that
\begin{align*}
R(\eta f)(x) -\eta(x) (R-G)f(x) = &\int (R(x,y)-g(x,y))(\eta(y)- \eta(x))f(y)d\mu(y) \\
&+ \int g(x,y)\eta(y)f(y)d\mu(y)\\
= &\int W_1(x,y) d\mu(y)+\int W_2(x,y) d\mu(y).
\end{align*}
Applying \eqref{integr} we easily estimate the summand that contains $W_2$. The function $W_1(x,y)$ vanishes if $|x-y|> 2\rho(y)$ or  $x,y \in I^c$. Therefore it can be verified that $W_1(x,y)=0$, if $x \notin 4I$ or $y \notin 4I$. Thus Lemma \ref{comut} follows by
\begin{align*}
\int_{4I} \Big| \int_{4I}  W_1(x,y)d\mu(y)\Big|d\mu(x) &\leq C \int_{4I} |f(y)|\( \int_{4I}\Big| \frac{1}{x^{\a+1}-y^{\a+1}}\Big| \frac{|x-y|}{\rho(z)} d\mu(x)
\) d\mu(y)\cr
&\leq C \int_{4I} |f(y)|d\mu(y).
\end{align*}
\end{proof}

\section{Proof of Theorem \ref{MainThm}}\label{Mainproof}
Before proving the main theorem we state a crucial consequence of Propositions \ref{lem2.2} and \ref{lem3.3}.
\begin{lema}\label{rty}
For $y_0>0$ we have
\begin{equation}\label{difference}
\sup_{y\in B\(y_0,\rho(y_0)\)} \ii \Big|R(x,y) - \wt r^{\rho(y_0)}(x,y) \Big| d\mu(x) \leq C,
\end{equation}
\end{lema}
\begin{proof}
By \eqref{prop_h} and \eqref{integr} we only need to establish that
$$\sup_{y\in B(y_0,\rho(y_0))} \ii \Big| \phi\(\frac{x-y}{\rho(y)}\)-\phi\(\frac{x-y}{\rho(y_0)}\) \Big| \ \Big|\frac{B}{x^{\a+1}-y^{a+1}}+ \frac{A-B}{x^{\a+1}+y^{a+1}} \Big| d\mu(x) \leq C.$$
In fact we will prove a stronger estimate, namely,
\begin{equation}\label{XYZ}
\sup_{y\in B(y_0,\rho(y_0))} \ii \Big| \phi\(\frac{x-y}{\rho(y)}\)-\phi\(\frac{x-y}{\rho(y_0)}\) \Big| \cdot \frac{1}{|x^{\a+1}-y^{a+1}|} d\mu(x) \leq C.
\end{equation}
Consider the case $y>y_0$ (if $y<y_0$  we use the same type of arguments). The integrant in \eqref{XYZ} is non-zero only when $3/2 \, \rho(y)<|x-y|<2\rho(y_0)$. But always $\rho(y_0)<2\rho(y)$ if $y\in B(y_0, \rho(y_0))$. Now, one can check that
$$ \sup_{y>0} \int^\8_{0} \chi_{\{3/2 \, \rho(y)<|x-y|<4\rho(y)\}}  \frac{1}{|x^{\a+1}-y^{\a+1}|} d\mu(x) \leq C,$$
which implies \eqref{XYZ}.
\end{proof}

{\em Proof of Theorem \ref{MainThm}.} Assume $f \in H_{at}^1(X)$. The operator $R:\LL \to \SSS '(X)$ is continuous (see \eqref{cont2}), so the first implication will be proved if we have established that there exists $C>0$ such that
$$\|R a\|_{\LL} \leq C$$
for any $H^1(X)-$atom $a$. Suppose $a$ is associated with $I = B(y_0, r)$ (recall that $r\leq \rho(y_0)$). We have
$$Ra = (Ra - \wt{r}^{\rho(y_0)}a) + \wt{r}^{\rho(y_0)}a.$$
The $\LL$-norm of the function $\wt{r}^{\rho(y_0)}a$ is bounded by a constant independent of $a$, because $a$ is also an $h^{1,\rho(y_0)}(X)$-atom (see Theorem \ref{localdecomp}). Therefore, the first part of the proof is finished by \eqref{difference}.

To prove the converse assume that $f, Rf\in \LL$. Introduce a family of intervals $\mc{I}= \{I_n=B(z_n,\rho(z_n))\}_{n=1}^{\8}$ such that $X = \bigcup_{n=1}^\8 I_n$ and $\mc I^* = \{4I: I\in \mc I\}$ has bounded overlap. Denote by $\eta_n$ a smooth partition of unity associated with the family $\mc{I}$, i.e.
$$\eta_n \in C^\8 \(0,\8\), \ \ \supp \ \eta_n \subset I_n, \ \ 0 \leq \eta_n \leq 1, \ \ \sum_{n=1}^\8 \eta_n (x) = \chi_{(0,\8)}(x), \ \ |\eta_n'(x)| \leq C\rho(z_n)^{-1}.$$
We are going to prove an atomic decomposition of $f=\sum_{n=1}^\8 \eta_n f$. Note that
$$\wt r^{\rho(z_n)} (f \eta_n) = \((\wt r^{\rho (z_n)}-R) (f \eta_n)\) + \(R(f \eta_n)-\eta_n (R-G)(f)\) - \eta_n \cdot G(f)+ \eta_n\cdot R(f).$$
By using \eqref{difference}, Lemma \ref{comut}, and \eqref{integr} we get
\begin{equation}\label{ooo}
\begin{split}
\sum_{n=1}^\8 \|\wt r^{\rho (z_n)} (f \eta_n)\|_{\LL}&\leq C\sum_{n=1}^{\8} \(\| \eta_n f\|_{\LL}+\| \chi_{4I_n} f\|_{\LL}+\|\eta_n Gf \|_{\LL}+\|\eta_n Rf \|_{\LL}\) \cr
&\leq C \(\|f\|_{\LL} +\|Rf\|_{\LL}\).
\end{split}
\end{equation}
Applying Theorem \ref{localdecomp}, we arrive at
\begin{equation}\label{yyy}
\eta_n \cdot f = \sum_{j=1}^\8 \lambda_{n,j} a_{n,j}, \text{ where }\sum_{j=1}^\8 |\lambda_{n,j}| \leq \|\wt r^{\rho (z_n)} (f\cdot \eta_n)\|_{\LL}, \text{ and } a_{n,j} \text{ are }h^{1,\rho(z_n)}(X)-\text{atoms}.
\end{equation}
From \eqref{ooo} and \eqref{yyy} we have obtained
\begin{equation}\label{dec3}
f=\sum_{n,j=1}^\8 \lambda_{n,j} a_{n,j} \ \text{ with } \quad \sum_{n,j=1}^\8 |\lambda_{n,j}| \leq C\(\|f\|_{\LL} +\|Rf\|_{\LL}\).
\end{equation}
Remark \ref{remark} states that $\supp \, a_{n,j} \subseteq 3 I_n$ for $j\geq 0$. Notice that for $y \in 3 \ I_n$ there exists $C>0$ such that

\begin{equation}\label{erttt}
\rho(z_n)/C \leq \rho(y) \leq C\rho(z_n) \quad \text { for all } n\geq 1 \text{ and } y \in I_n.
\end{equation}

Because of this, each $a_{n,j}$ can be decomposed into a sum of at most $k$ atoms of the space $H^1(X)$ (where the number $k$ depends only on $\a$ and the constant $C$ from \eqref{erttt}). Finally, Theorem \ref{MainThm} follows by applying \eqref{dec3}.
\begin{flushright}
$\square$
\end{flushright}

\section{Auxiliary estimates}\label{aux}
This section is devoted to proving Propositions \ref{lem2.2} and \ref{lem3.3}. The letters $c,C,N,M$ will denote positive constants ($N,M$ are arbitrarily large). We also use the convention that $\int_p^q \dots=0$, when $p \geq q$. For further references we figure out some properties of the Bessel function $I_\nu \ (\nu >0)$ (see, e.g., \cite{Watson}):
\begin{align}
\label{derrI}
&\dx \( x^{-\nu}I_{\nu}(x)\) =x^{-\nu} I_{\nu+1}(x) &&\text{ for } x>0 ,\\
\label{zeroI}
&0 < I_\nu (x) = 2^{-\nu} \Gamma(\nu+1)^{-1} x^\nu + O(x^{\nu + 2}) &&\text{ for }0<x<C,\\
\label{infI}
&U_\nu (x) = (2\pi)^{-1\sl 2} + O(x^{-1}) &&\text{ for } x>C,
\end{align}
where
$$U_\nu (x)=I_\nu (x) e^{-x} \sqrt{x}.$$

\subsection{Proof of Proposition \ref{lem2.2}}\label{proof2.2}
\begin{proof}
Assume $y=1$. By using \eqref{kerbess} and \eqref{derrI} we get
\begin{align} \label{kernel1}
\wt R (x,1) =&  - \ii (2t)^{-2} \exp\left(-\frac{x^2+1}{4t}\right)
x^{-(\al-3)/2}I_{(\al-1)/2} \left(\frac{x}{2t}\right) \dtt \cr &+
\ii (2t)^{-2} \exp\left(-\frac{x^2+1}{4t}\right)
x^{-(\al-1)/2}I_{(\al+1)/2} \left(\frac{x}{2t}\right) \dtt \cr = &
\ii \wt Q_1(x,t) dt + \ii \wt Q_2(x,t) dt.
\end{align}

In calculations below we will often use the following formula:
$$\exp\({-\frac{x^2+1}{4t}}\)I_\nu\(\frac{x}{2t}\)\sqrt{\frac{x}{2t}}=\exp\({-\frac{(x-1)^2}{4t}}\)U_\nu\(\frac{x}{2t}\).$$

Define
$$h(x) = \wt R (x,1) - \frac{A-B}{x^{\a+1}+1}- \frac{B}{x^{\a+1}-1}.$$

To prove \eqref{prop_h} we consider three cases.

{\bf Case 1:} $x> 3/2$.\\
Under this assumption $x-1 \sim x$. Then we get estimates:
\begin{align*}
\int_0^x |\wt Q_2(x,t)| dt &\leq C \,\int_0^x t^{-2} \exp\( -\frac{(x-1)^2}{4t}\) U_{(\al+1)\slash 2 } \(\xt\)x^{-\al\slash 2} dt\\
&\leq C\int_0^x t^{-2} \(\frac{t}{x^2}\)^N x^{-\al\slash 2} dt \leq C x^{-M},\\
\int_x^{x^2} |\wt Q_2(x,t)| dt &\leq C \,\int_x^{x^2} t^{-2} \(\frac{t}{x^2}\)^N x^{(1-\al)\slash 2}\(\frac{x}{t}\)^{\frac{\al+1}{2}} \dtt \leq C \int_0^{x^2} \frac{t ^{N-3-\al\slash 2}}{x^{2N-1}} dt \leq C x^{-\al-3},\\
\int_{x^2}^\8 |\wt Q_2(x,t)| dt &\leq C \,\int_{x^2}^\8 t^{-2} x^{-(\al-1)\slash 2}\(\frac{x}{t}\)^{(\al+1)\slash 2} \dtt \leq C \int_{x^2}^\8 x t^{-3-\al\slash 2} dt \leq C x^{-\al-3},
\end{align*}
which imply
\begin{equation}\label{infi1}
\ii| \wt Q_2(x,t)| dt \leq C \, x ^{-\al-3}.
\end{equation}

Our next task is to obtain
\begin{equation}\label{infinity}
\Big| \ii \wt Q_1(x,t) dt - \frac{A-B}{x^{\al+1}+1} -\frac{B}{x^{\al+1}-1}  \Big| \leq C x^{-\al-2}.
\end{equation}

By using the same methods as we have utilized to estimate the integral $\int_0^x |\wt Q_2(x,t)| dt$ we deduce
\begin{equation}\label{ppt}
\int_0^x| \wt Q_1(x,t)| dt \leq C x^{-M}.
\end{equation}
Moreover,
\begin{align}
\Big| \int_x^{\8} \wt Q_1(x,t) dt &- \frac{A}{x^{\al+1}}\Big|
=\Big|-\int_x^\8 (2t)^{-2} \exp\(-\frac{x^2+1}{4t}\) x^{-\frac{\a-3}{2}}I_{(\a-1)/2} \left(\frac{x}{2t}\right) \dtt\cr
&+\ii (2t)^{-2} \exp\(-\frac{x^2}{4t}\) x^{-\frac{\a-3}{2}} \frac{2^{-(\a-1)/2}}{\Gamma\(\frac{\a+1}{2}\)}\left(\frac{x}{2t}\right)^{\frac{\a-1}{2}} \dtt\Big| \cr
\leq &\Big|\int_x^\8 (2t)^{-2} \(\exp\(-\frac{x^2+1}{4t}\)- \exp\(-\frac{x^2}{4t}\)\)x^{-\frac{\a-3}{2}}I_{(\a-1)/2} \left(\frac{x}{2t}\right) \dtt\Big|\label{est133} \\
&+\Big|\int_x^\8 (2t)^{-2} \exp\(-\frac{x^2}{4t}\)x^{-\frac{\a-3}{2}}\(I_{(\a-1)/2} \left(\frac{x}{2t}\right) -\frac{2^{-(\a-1)/2}}{\Gamma\(\frac{\a+1}{2}\)}\left(\frac{x}{2t}\right)^{\frac{\a-1}{2}} \) \dtt\Big|\cr
&+\Big|\int_0^x (2t)^{-2} \exp\(-\frac{x^2}{4t}\) x^{-\frac{\a-3}{2}} \frac{2^{-(\a-1)/2}}{\Gamma\(\frac{\a+1}{2}\)}\left(\frac{x}{2t}\right)^{\frac{\a-1}{2}} \dtt\Big|\notag
\end{align}
and
\begin{equation}\label{rrr}
\Big| \frac{A}{x^{\al+1}} - \frac{A-B}{x^{\al+1}+1}-\frac{B}{x^{\al+1}-1}\Big| \leq C x^{-2\al-2} \qquad \text{ for } x>3/2.
\end{equation}
Applying \eqref{zeroI} and the mean-value theorem to \eqref{est133}, we get
\begin{equation}\label{nnnhhh}
\Big| \int_x^{\8} \wt Q_1(x,t) dt - \frac{A}{x^{\al+1}}\Big| \leq C x^{-\al-3}.
\end{equation}
Now \eqref{infinity} is a consequence of \eqref{ppt}, \eqref{rrr}, and \eqref{nnnhhh}. From \eqref{infi1}--\eqref{infinity} we conclude that
\begin{equation}\label{summ1}
\int_{3/2}^\8 |h(x)|d\mu(x) = \int_{3/2}^\8 \Big|\wt R(x,1) -  \frac{A-B}{x^{\al+1}+1}-\frac{B}{x^{\al+1}-1}\Big| d\mu(x) \leq C.
\end{equation}

{\bf Case 2:} $x< 1/2$.\\
From \eqref{zeroI}--\eqref{infI} it follows:
\begin{align*}
&\int_0^x |\wt Q_{1}(x,t)| dt \leq C\int_0^x \frac{x^{1-\al\slash 2}}{t^{2}} \exp\left(-\frac{(x-1)^2}{4t}\right) U_{\frac{\al -1}{2}}\(\xt\) dt \leq C x^{1-\al\sl 2} \int_0^x t^{N-2} dt \leq Cx,\\
&\int_x^1 |\wt Q_{1}(x,t)| dt \leq C\int_x^1 t^{-2} \exp\left(-\frac{x^2+1}{4t}\right) x^{-\frac{\al-3}{2}} \(\frac{x}{t}\)^{(\al-1)\sl 2} \dtt \leq C x \int_0^1 t^{M} dt \leq Cx,\\
&\int_1^\8 |\wt Q_{1}(x,t)| dt \leq C\int_1^\8 t^{-2} \exp\left(-\frac{x^2+1}{4t}\right) x^{\frac{3-\al}{2}} \(\frac{x}{t}\)^{\frac{\al-1}{2}} \dtt
\leq C x \int_1^\8 t^{-2-\al \sl 2} dt \leq Cx,
\end{align*}

Thus $\ii |\wt Q_1(x,t)|dt \leq Cx$. By the same arguments we also obtain
$\ii |\wt Q_2(x,t)|dt \leq Cx$. Hence, $|\wt R (x,1)|\leq Cx$. As a consequence, for $x<1/2$, we have
\begin{eqnarray}\label{lab1}
|h(x)+A-2B| = \Big|\wt R (x,1)-\frac{A-B}{x^{\al+1}+1}-\frac{B}{x^{\al+1}-1}+A-2B\Big| \leq Cx,\\ \label{summ2}
\int_{0}^{1/2} |h(x)|d\mu(x) \leq \int_{0}^{1/2} \(|\wt R(x,1)| + \Big| \frac{A-B}{x^{\al+1}+1}+\frac{B}{x^{\al+1}-1}\Big| \) d\mu(x) \leq C.
\end{eqnarray}

{\bf Case 3:} $1/2<x<3/2$.\\
In this case a slightly different form of \eqref{kernel1} is needed, i.e.,
\begin{align}
\ii \dx \wt T_t (x,1) \dtt &= -(x-1) \ii (2t)^{-2} \exp\(-\frac{x^2+1}{4t}\) x^{-(\al-1)/2}I_{(\al-1)/2} \(\frac{x}{2t}\) \dtt\notag \\
+ &\ii (2t)^{-2} \exp\(-\frac{x^2+1}{4t}\) x^{-(\al-1)/2} \( I_{(\al+1)/2}\(\frac{x}{2t}\)- I_{(\al-1)/2} \(\frac{x}{2t}\)\) \dtt  \cr
&= \ii \wt Q_3(x,t) dt + \ii \wt Q_4(x,t) dt. \label{kernel2}
\end{align}

We claim that
\begin{equation}\label{eee}
\ii |\wt Q_4(x,t)| dt\leq C|x-1|^{-1\sl 2}.
\end{equation}
Indeed, by using \eqref{zeroI} and \eqref{infI} we get
\begin{align*}
\int_0^1 \Big| \wt Q_4(x,t)\Big| dt &\leq C \int_0^1 t^{-2}
\exp\(-\frac{(x-1)^2}{4t}\) x^{-\al\sl 2} \Big| U_{(\al+1)\sl 2}
\(\xt\) - U_{(\al-1)\sl 2} \(\xt\)\Big|dt \cr &\leq C \int_0^1
t^{-2} \left( \frac{t}{(x-1)^2}\right)^{1\sl 4} \frac{t}{x} dt \leq
C |x-1|^{-1\sl 2},\cr
\int_1^\8\Big| \wt Q_4(x,t)\Big| dt &\leq C \int_1^\8 t^{-2}
x^{\frac{1-\al}{2}} \( \(\frac{x}{t}\)^{\frac{\al+1}{2}} +
\(\frac{x}{t}\)^{(\al-1)\sl 2} \)\dtt \leq C \int_1^\8 t^{-2-\al\sl
2} dt\leq C.
\end{align*}

Next, observe that
\begin{equation}\label{fff}
\int_1^\8 \Big| \wt Q_3(x,t) \Big|dt \leq C|x-1| \int_1^\8 t^{-2} x^{-(\al-1)\sl 2} \(\frac{x}{t}\)^{(\al-1)\sl 2} \dtt \leq C.
\end{equation}

Moreover,
\begin{align}
\Big| \int_0^1 \wt Q_3(x,t) dt - \frac{B(\a+1)^{-1}}{x^{\al/2}(x-1)}\Big| \leq \Big|\int_1^\8 \frac{\sqrt{2}(x-1)}{4t}\exp\(-\frac{(x-1)^2}{4t}\) x^{-\a /2} \frac{1}{\sqrt{2\pi }} \dttt\Big| \label{gggg} \\
+\Big| \int_0^1 \frac{\sqrt{2}(x-1)}{4t} \exp\(-\frac{(x-1)^2}{4t}\) x^{-\a/2} \(U_{(\al-1)\sl 2}\(\xt\)-\frac{1}{\sqrt{2\pi}}\)  \dttt\Big|.\notag
\end{align}
Applying \eqref{infI} to \eqref{gggg} we deduce
\begin{equation}\label{rrss}
\Big| \int_0^1 \wt Q_3(x,t) dt - \frac{B(\a+1)^{-1}}{x^{\al/2}(x-1)}\Big| \leq C.
\end{equation}
One can easily check that
\begin{align}\label{rrttt}
\Big| \frac{B(\a+1)^{-1}}{x^{\al/2}(x-1)} - \frac{B}{x^{\a+1}-1}-\frac{A-B}{x^{\a+1}+1}\Big| \leq C.
\end{align}

From \eqref{eee}, \eqref{fff}, \eqref{rrss}, and \eqref{rrttt} we conclude
\begin{equation}\label{summ3}
\int_{1/2}^{3/2} |h(x)|d\mu(x) = \int_{1/2}^{3/2} \Big|\wt R (x,1) - \frac{B}{x^{\a+1}-1}-\frac{A-B}{x^{\a+1}+1}\Big| d\mu(x) \leq C.
\end{equation}

Finally, as a consequence of \eqref{summ1}--\eqref{summ2}, and \eqref{summ3} we obtain that $h$ satisfies desired properties \eqref{prop_h}.  The proposition in the general case of $y >0$ follows by applying the homogeneity
\begin{equation}\label{dill}
\wt R(x,y) = y^{-\a-1}\wt R\(\frac{x}{y}, 1\).
\end{equation}
\end{proof}
\subsection{Proof of Proposition \ref{lem3.3}}\label{proof3.3}

\begin{proof}
Let us set
\begin{equation}\label{gdef}
g(x,y) = R(x,y) - \phi\(\frac{x-y}{\rho(y)}\) \( \frac{B}{x^{\a+1}-y^{\a+1}} - \frac{A-B}{x^{\a+1}+y^{\a+1}} \).
\end{equation}
We will prove that \eqref{integr} is satisfied. By using \eqref{semi} and \eqref{derrI} we get
\begin{equation}\label{Rieszdecomp}
R(x,y) = \ii T_t^{[1]}(x,y) \dtt + \ii T_t^{[2]}(x,y) \dtt =\ii T_t^{[3]}(x,y) \dtt +\ii T_t^{[4]}(x,y) \dtt,
\end{equation}
where
\begin{eqnarray*}
&T_t^{[1]}(x,y) =& \(\frac{2e^{-2t}}{1-e^{-4t}}\)^2 y (xy)^{-\frac{\a-1}{2}} \exp\(-\frac{1+e^{-4t}}{1-e^{-4t}} \frac{x^2+y^2}{2}\) I_{\frac{\a+1}{2}}\(\frac{2e^{-2t}}{1-e^{-4t}}xy\),\cr
&T_t^{[2]}(x,y) =&- \frac{2e^{-2t}(1+e^{-4t})}{(1-e^{-4t})^2} x (xy)^{-\frac{\a-1}{2}} \exp\(-\frac{1+e^{-4t}}{1-e^{-4t}} \frac{x^2+y^2}{2}\) I_{\frac{\a-1}{2}}\(\frac{2e^{-2t}}{1-e^{-4t}}xy\),\cr
&T_t^{[3]}(x,y) =&- \frac{2e^{-2t}(1+e^{-4t})}{(1-e^{-4t})^2} (xy)^{-\frac{\a-1}{2}}(x-y) \exp\(-\frac{1+e^{-4t}}{1-e^{-4t}} \frac{x^2+y^2}{2}\) I_{\frac{\a-1}{2}}\(\frac{2e^{-2t}xy}{1-e^{-4t}}\),\cr
&T_t^{[4]}(x,y) =&\frac{2e^{-2t}}{1-e^{-4t}} y (xy)^{-\frac{\a-1}{2}} \exp\(-\frac{1+e^{-4t}}{1-e^{-4t}}  \frac{x^2+y^2}{2}\) \cr
&&\cdot\Bigg(\frac{2e^{-2t}}{1-e^{-4t}} I_{\frac{\a+1}{2}}\(\frac{2e^{-2t}}{1-e^{-4t}}xy\)
-\frac{1+e^{-4t}}{1-e^{-4t}} I_{\frac{\a-1}{2}}\(\frac{2e^{-2t}}{1-e^{-4t}}xy\)\Bigg).
\end{eqnarray*}

Note that
\begin{align}
&\exp\(-\frac{1+e^{-4t}}{1-e^{-4t}} \frac{x^2+y^2}{2}\) I_{\mu}\(\frac{2e^{-2t}}{1-e^{-4t}}xy\)\(\frac{2e^{-2t}xy}{1-e^{-4t}}\)^{1\sl 2}
\label{IchangeU} \\
&= \exp\(-\frac{1+e^{-4t}}{1-e^{-4t}} \frac{(x-y)^2}{2}\) \exp\(-\frac{(1-e^{-2t})^2}{1-e^{-4t}}xy\)U_{\mu}\(\frac{2e^{-2t}xy}{1-e^{-4t}}\).\notag
\end{align}

The formula \eqref{IchangeU} will be frequently used, without additional comments, when we deal with $I_\mu(\theta)$ for $\theta >C$.

We provide the proof in six cases as it is shown in Figure 1 on page \pageref{page_jpg}.
\begin{figure}[h]\label{page_jpg}
\includegraphics[height=3in]{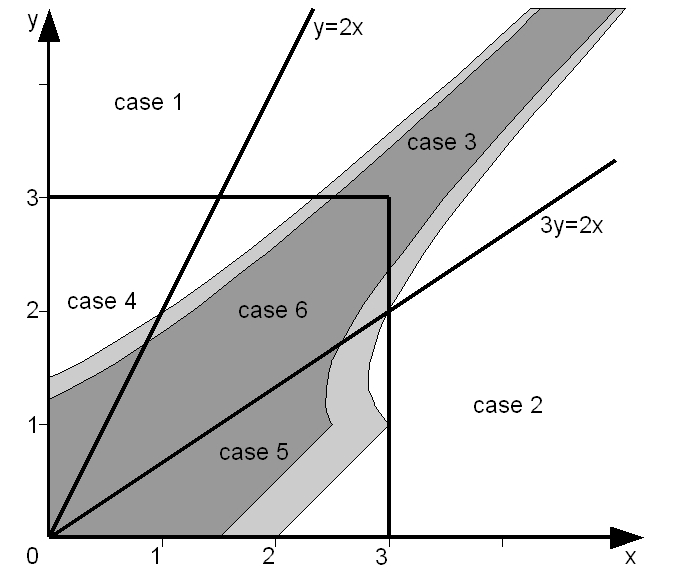}\\
\caption{The partition of $X \times X$}
\end{figure}
The grey part denotes the support of $\phi\((x-y)/ \rho(y)\)$. Moreover, the dark grey color means that $\phi\((x-y) /\rho(y)\)=1$.

In Cases 1, 2, 4, 5 we will use the decomposition \eqref{Rieszdecomp} that contains $T^{[1]}_t$ and $T_t^{[2]}$.

{\bf Case 1: $y>3$, $x<y\sl 2$.}\\
At the beginning we consider $T_t^{[1]}$  and $t<1$. Under additional assumption $xy<1$ we get
\begin{eqnarray*}
&\int_{xy}^{1}|T_t^{[1]}(x,y)| \dtt &\leq
C \int_{xy}^{1} t^{-2}y(xy)^{-\frac{\a-1}{2}}
\exp\(- \frac{y^2}{ct}\)
\(\frac{xy}{t}\)^{\frac{\a+1}{2}}\dtt
\leq C y^{-M},\cr
&\int_0^{xy}|T_t^{[1]}(x,y)| \dtt &\leq C \int_0^{xy}
\frac{y}{t^2}(xy)^{-\frac{\a-1}{2}} \exp\(-\frac{y^2}{ct}\) \exp\(-ctxy\) U_{\frac{\a+1}{2}}\(\frac{2e^{-2t}xy}{1-e^{-4t}}\) \sqrt{\frac{t}{xy}}\dtt \cr
&&  \leq C \int_0^{xy} t^{-2} y(xy)^{-\frac{\a}{2}}\(\frac{t}{y^2}\)^{N} dt \leq C\( \frac{x}{y}\)^N x^{-\a-1}.\label{est1}
\end{eqnarray*}
In the last line we have used \eqref{IchangeU} and \eqref{infI}. If $xy>1$ we similarly get $\int_0^1 |T_t^{[1]}(x,y)| \dtt \leq C y^{-M}$.

Next, we deal with $T_t^{[1]}$ and $t>1$. If $xy>e^2$ then
\begin{equation*}
\begin{split}
\int_1^{\log \sqrt{xy}} |T_t^{[1]}(x,y)| \dtt \leq C
\int_1^{\log \sqrt{xy}}e^{-4t} y (xy)^{-(\a-1)\sl 2} \exp (-cy^2)
(e^{-2t} xy)^{-1\sl 2} \dtt \leq  C y^{-M}, \\
\int_{ \log \sqrt{xy}}^{\8}|T_t^{[1]}(x,y)| \dtt \leq C \int_{\log
\sqrt{xy}}^{\8} e^{-4t} y (xy)^{-(\a-1)\sl 2} \exp(-cy^2) (e^{-2t}
xy)^{(\a+1)\sl 2} \dtt \leq C  y^{-M}.
\end{split}
\end{equation*}
Identically, when $xy<e^2$ we have $\int_1^\8 |T_t^{[1]}(x,y)| \dtt \leq Cy^{-M}$.

We can write the same estimates for $T_t^{[2]}$. Thus we get
$$\ii \(|T_t^{[1]}(x,y)|+|T_t^{[2]}(x,y)|\)\dtt \leq C y^{-M} \max (1,x^{M-\al-1}).$$
Observe that $g(x,y)=R(x,y)$ (see Figure 1), so the last estimate implies
\begin{equation}\label{all1}
\sup_{y>3} \int_0^{y\sl 2} |g(x,y)| d\mu(x) < \8.
\end{equation}

{\bf Case 2: $x>3$, $y<2x\sl 3$.}\\
We proceed very similarly to Case $1$ and obtain
\begin{equation*}
\int_0^{\8}\(|T_t^{[1]}(x,y)|+|T_t^{[2]}(x,y)|\) \dtt \leq C x^{-M} \max ( 1, y^{M-\a-1} ).
\end{equation*}
We have $g(x,y)=R(x,y)$ (see Figure 1). Hence
\begin{equation}\label{all2}
\sup_{y<2} \int_{3}^{\8} |g(x,y)| d\mu(x) < \8, \qquad
\sup_{y\geq 2} \int_{3y/2}^{\8} |g(x,y)| d\mu(x) < \8.
\end{equation}

{\bf Case 3:} ($x>3$ or $y>3$) and $|x-y| < y\sl 2$.\\
Notice that
\begin{align}
T_t^{[4]} (x,y) = &\(\frac{2e^{-2t}}{1-e^{-4t}}\)^{1\sl 2} y (xy)^{-\a\sl 2} \exp \(-\frac{1+e^{-4t}}{1-e^{-4t}} \frac{(x-y)^2}{2} \) \label{decompV} \\
&\times \exp\(-\frac{(1-e^{-2t})^2}{1-e^{-4t}}xy\)  \(V_t^{[4']} (x,y) +V_t^{[4'']} (x,y) +V_t^{[4''']} (x,y)\),\notag
\end{align}
where
\begin{eqnarray*}
&V_t^{[4']} (x,y) &=\frac{2e^{-2t}}{1-e^{-4t}} \(U_{(\a+1)\sl 2} \(\frac{2e^{-2t}}{1-e^{-4t}}xy\) - \frac{1}{\sqrt{2 \pi}}\),\cr
&V_t^{[4'']} (x,y) &= \frac{1}{\sqrt{2 \pi}}\(\frac{2e^{-2t}}{1-e^{-4t}} - \frac{1+e^{-4t}}{1-e^{-4t}}\)=-\frac{(1-e^{-2t})^2}{\sqrt{2\pi }(1-e^{-4t})},\cr
&V_t^{[4''']} (x,y) &= -\frac{1+e^{-4t}}{1-e^{-4t}} \(U_{(\a-1)\sl 2} \(\frac{2e^{-2t}}{1-e^{-4t}}xy\) - \frac{1}{\sqrt{2 \pi}}\).
\end{eqnarray*}

By using \eqref{decompV} and \eqref{infI} one obtains
\begin{equation} \label{est2}
\int_0^1 |T_t^{[4]}(x,y)|\dtt \leq C |x-y|^{-1\sl 2} x^{-\a-1}.
\end{equation}

Also, as in Case $1$, we get
\begin{equation} \label {est3}
\int_1^{\8} \(|T_t^{[3]}(x,y)| + |T_t^{[4]}(x,y)|\) \dtt \leq C x^{-M}.
\end{equation}

Next,
\begin{equation} \label{decompD}
\int_0^{1} T_t^{[3]}(x,y)\dtt - \chi_{\{y|x-y|<1\}} \frac{B(xy)^{-\a/2}}{(\al+1) (x-y)} = D_1 - D_6=\sum_{j=1}^{5} (D_j - D_{j+1}),
\end{equation}
where
\begin{eqnarray*}
D_2 &=& -\int_0^1 \frac{e^{-t}(1+e^{-4t})}{\sqrt{\pi}(1-e^{-4t})^{3\sl 2}} \frac{x-y}{(xy)^{\frac{\a}{2}}} \exp\(-\frac{1+e^{-4t}}{1-e^{-4t}} \frac{(x-y)^2}{2}\)
\exp\(-\frac{(1-e^{-2t})^2}{1-e^{-4t}}xy\) \dtt , \cr
D_3 &=& -\int_0^1 \frac{1}{4t \sqrt{\pi} } (xy)^{-\a\sl 2}(x-y) \exp\(- \frac{(x-y)^2}{4t}\)
\exp\(-txy\) \dttt,\cr
D_4 &=& -\int_0^{y^{-2}/4} \frac{1}{4t \sqrt{\pi} } (xy)^{-\a\sl 2}(x-y) \exp\(- \frac{(x-y)^2}{4t}\)
\exp\(-txy\) \dttt,\cr
D_5 &=& -\int_0^{y^{-2}/4} \frac{1}{4t \sqrt{\pi}} (xy)^{-\a\sl 2}(x-y) \exp\(- \frac{(x-y)^2}{4t}\) \dttt,\cr
D_6 &=& -\chi_{\{y|x-y|<1\}}\int_0^{\8} \frac{1}{4t \sqrt{\pi}} \frac{x-y}{(xy)^{\frac{\a}{2}}} \exp\(- \frac{(x-y)^2}{4t}\) \dttt= -\chi_{\{y|x-y|<1\}} \frac{(xy)^{-\a/2}}{\sqrt{\pi} (x-y)}.
\end{eqnarray*}

By using the mean-value theorem, \eqref{infI}, and \eqref{IchangeU} one obtains
\begin{equation}\label{est5}
|D_j-D_{j+1}| \leq C x^{-\a-1/2} |x-y|^{-1/2} \quad \text{ for } j=1,2.
\end{equation}

To deal with $D_j-D_{j+1}$ for $j=3,4,5$ we consider:

{\bf Subcase 1: $y|x-y|<1$.}
\begin{align}
&|D_3-D_4| \leq C \int_{y^{-2}/4}^1 t^{-1} x^{-\a} |x-y| \dttt \leq Cx^{-\a+2}|x-y|, \notag \\
&|D_4-D_5| = \int_0^{y^{-2}/4} = \int_0^{(x-y)^{2}/4} ... \ dt +\int_{(x-y)^2/4}^{y^{-2}/4} ...  \ dt = Y_1 + Y_2,\cr
&|Y_1| \leq C \int_0^{(x-y)^{2}} t^{N-1} x^{-\a+2} |x-y|^{1-2N}  dt\leq C x^{-\a+2}|x-y|, \label{est6} \\
&|Y_2| \leq C \int_{(x-y)^2/4}^{y^{-2}/4} x^{-\a+2} |x-y| \dttt\leq Cx^{-\a+2} |x-y| \ln \frac{1}{y|x-y|},\cr
&|D_5-D_6| \leq C \int_{y^{-2}/4}^\8 t^{-1} x^{-\a} |x-y| \dttt \leq Cx^{-\a+2}|x-y|.\notag
\end{align}

{\bf Subcase 2: $y|x-y|>1$.}
\begin{align}
&|D_3-D_4| = \int_{y^{-2}/4}^1 = \int_{y^{-2}/4}^{y^{-1}|x-y|/4}... \ dt +\int_{y^{-1}|x-y|/4}^1 ... \ dt = Y_3 + Y_4,\cr
&|Y_3| \leq C \int_{y^{-2}/4}^{y^{-1}|x-y|/4} t^{N-2} x^{-\a} |x-y|^{1-2N} dt\leq C x^{-\a+1-N}|x-y|^{-N},\cr
&|Y_4| \leq C \int_{y^{-1}|x-y|/4}^{1} t^{-2} x^{-\a}(txy)^{-N} |x-y| dt\leq Cx^{-\a+1-N}|x-y|^{-N}, \label{est7}\\
&|D_4-D_5| \leq C \int_{0}^{y^{-2}/4} t^{N-1} x^{-\a+2} |x-y|^{1-2N} dt \leq Cx^{-\a+1-M} |x-y|^{-M}, \cr
&|D_5-D_6|=|D_5| = C \int_{0}^{y^{-2}/4} t^{N-1} x^{-\a} |x-y|^{1-2N} \dttt \leq Cx^{-\a+1-M} |x-y|^{-M}.\notag
\end{align}

Reassuming, \eqref{decompV}--\eqref{est7} lead to
\begin{equation}\label{all33}
\sup_{y>2} \int_{y\sl 2}^{3y/2} \Big| R(x,y) -\chi_{\{y|x-y|<1\}}(x) \frac{B(xy)^{-\a/2}}{(\al+1) (x-y)}\Big|  d\mu(x) < \8.
\end{equation}
Moreover,
\begin{equation}\label{edrt}
\chi_{\{y|x-y|<2\}}(x)\Big| \frac{B(xy)^{-\a/2}}{(\al+1) (x-y)} -\frac{B}{x^{\al+1}-y^{\al+1}}-\frac{A-B}{x^{\al+1}+y^{\al+1}}\Big| \leq C x^{-\al-1}.
\end{equation}
We claim that
\begin{equation}\label{all3}
\sup_{2<y<3}\int_{3}^{3y/2}|g(x,y)|d\mu(x) \leq C \quad \text{and} \quad \sup_{3<y}\int_{y/2}^{3y/2}|g(x,y)|d\mu(x) \leq C.
\end{equation}
To prove \eqref{all3} we split the area of integration into three parts that correspond to white, light grey, and dark grey regions from Figure 1.
\begin{itemize}
\item if $y|x-y|>2$ we have $\phi\((x-y)/\rho(y)\)=0$ and we deduce the statement directly from \eqref{all33}.
\item if $1\leq y|x-y|\leq 2$ then we apply \eqref{all33}--\eqref{edrt}, and the inequality
$$\sup_{y>2}\int_{1<y|x-y|<2}\(\Big|\frac{B}{x^{\al+1}-y^{\al+1}}\Big|+\Big|\frac{A-B}{x^{\al+1}+y^{\al+1}}\Big|\)d\mu(x)\leq C.$$
\item if $y|x-y|<1$ then $\phi\((x-y)/\rho(y)\)=1$  and we use again \eqref{all33}--\eqref{edrt}.
\end{itemize}

{\bf Case 4:} $x,y<3$, $x <y\sl 2$.\\
By similar analysis to that we have used in Case 1 we obtain
\begin{align*}
&\int_{0}^{xy} \(|T_t^{[1]} (x,y)|+| T_t^{[2]} (x,y)|\) \dtt \leq C \(\frac{x}{y}\)^{M} x^{-\a-1}, \cr
&\int_{xy}^{1} \(|T_t^{[1]} (x,y)|+ |T_t^{[2]} (x,y)|\) \dtt \leq C (xy)^{-\a\sl 2} x^{-1},\cr
&\int_{1}^{\8} \(|T_t^{[1]} (x,y)|+ |T_t^{[2]} (x,y)|\) \dtt \leq C.
\end{align*}

Therefore
\begin{equation*}
\sup_{y<3} \int_0^{y\sl 2} |R(x,y)| d\mu(x) < \8.
\end{equation*}

and, consequently,
\begin{equation}\label{all4}
\sup_{y<3} \int_0^{y\sl 2} |g(x,y)| d\mu(x) < \8,
\end{equation}
since
$$ \sup_{y<3} \int_0^{y\sl 2}\( \Big| \frac{B}{x^{\a+1}-y^{\a+1}}\Big| + \Big|\frac{A-B}{x^{\a+1}+y^{\a+1}} \Big| \)\, d\mu(x) <\8  .$$

{\bf Case 5:} $x,y<3$, $y<2x/3$.\\
By using \eqref{zeroI} and \eqref{infI}, similarly as in Case 2, one obtains
\begin{align}
&\int_{0}^{xy} \(|T_t^{[1]} (x,y)|+|T_t^{[2]} (x,y)|\) \dtt \leq C \(\frac{y}{x}\)^{M} x^{-\a-1},\notag \\
&\int_{xy}^{x^2} |T_t^{[1]} (x,y)| \dtt \leq C \(\frac{y}{x}\)^{2} x^{-\a-1},\label{est11} \\
&\int_{x^2}^{1} |T_t^{[1]} (x,y)| \dtt \leq C \(\frac{y}{x}\)^{2} x^{-\a-1}, \notag \\
&\int_{1}^{\8} \(|T_t^{[1]}(x,y)| + |T_t^{[2]}(x,y)|\) \dtt \leq C.\notag
\end{align}

Recall that $A=-2 \gamma_1 \gamma_2^{-1}$, where $\gamma_1 = \Gamma(\a/2+1)$ and $\gamma_2 = \Gamma((\a+1)/2)$. We write

\begin{equation} \label{est12}
\int_{xy}^{1} T_t^{[2]}(x,y)\dtt -\frac{Ax}{(x^2+y^2)^{\a/2+1}} = E_1 - E_4= \sum_{j=1}^{3} (E_j - E_{j+1}),
\end{equation}
where
\begin{eqnarray*}
E_2 &=& -\int_{xy}^1 \frac{2e^{-2t}(1+e^{-4t})}{(1-e^{-4t})^{2}} x  \exp\(-\frac{1+e^{-4t}}{1-e^{-4t}} \frac{x^2+y^2}{2}\) \gamma_2^{-1} \(\frac{e^{-2t}}{1-e^{-4t}}\)^{\frac{\a-1}{2}}
 \dtt , \cr
E_3 &=& -\int_{xy}^1   \frac{2}{\gamma_2} (4t)^{-\a/2-1} x  \exp\(-\frac{x^2+y^2}{4t} \)
 \dttt , \cr
E_4 &=& -\int_{0}^\8 \frac{2}{\gamma_2} (4t)^{-\a/2-1}  x \exp\(-\frac{x^2+y^2}{4t}\)
 \dttt =-2\frac{\gamma_1}{\gamma_2} \frac{x}{(x^2+y^2)^{\a/2+1}}.
\end{eqnarray*}
Applying \eqref{zeroI} and the mean-value theorem, one gets
\begin{eqnarray} \label{est13}
&|E_1-E_2|&\leq C\, y^2 x^{-\a-3},\cr
&|E_2-E_{3}|&\leq C\, x^{-\a+1},\cr
&|E_3-E_{4}|&\leq C \max (1,y^M x^{-\a-1-M}).
\end{eqnarray}
Moreover,
\begin{equation}\label{est131}
\Big|E_4 - \frac{B}{x^{\al+1}-y^{\al+1}} +\frac{A-B}{x^{\al+1}+y^{\al+1}}\Big| \leq Cy x^{-\al-2}.
\end{equation}
As a consequence of \eqref{est11}--\eqref{est131} we get
\begin{equation}\label{pre1}
\sup_{y<2} \int_{3y/2}^{3}\Big|R(x,y)- \frac{B}{x^{\al+1}-y^{\al+1}} -\frac{A-B}{x^{\al+1}+y^{\al+1}}\Big| d\mu(x) < \8.
\end{equation}
Also,
\begin{equation}\label{pre2}
\sup_{y<2} \int_{1}^{3}\chi_{\{y<2x/3\}}\(\Big|\frac{B}{x^{\al+1}-y^{\al+1}} \Big| +\Big|\frac{A-B}{x^{\al+1}+y^{\al+1}}\Big|\) d\mu(x) < \8.
\end{equation}
Observe that if $x<1$ then $\phi((x-y)/\rho(y))=1$ (see Figure 1). Therefore \eqref{pre1}-\eqref{pre2} lead to
\begin{eqnarray}\label{all5}
\sup_{y<2} \int_{3y/2}^{3}|g(x,y)| d\mu(x) < \8.
\end{eqnarray}

{\bf Case 6:} $x,y<3$, $|x-y| < y/2$.\\
By using the decomposition
\eqref{decompV} one obtains
\begin{equation} \label{est9}
\int_0^{xy} |T_t^{[4]}(x,y)|\dtt \leq C |x-y|^{-1/2} x^{-\a-1/2}.
\end{equation}

In addition
\begin{align} \label{est10}
\int_{xy}^{1} \( |T_t^{[3]}(x,y)|+|T_t^{[4]}(x,y)|\) \dtt \leq C x^{-\a-1}, \qquad
\int_{1}^{\8} \( |T_t^{[3]}(x,y)|+|T_t^{[4]}(x,y)|\) \dtt \leq C.
\end{align}

Denote
$$\int_0^{xy} T_t^{[3]}(x,y)\dtt - \(\frac{B}{x^{\a+1}-y^{\a+1}}+ \frac{A-B}{x^{\a+1}+y^{\a+1}}\) = F_1 - F_5= \sum_{j=1}^{4} (F_j - F_{j+1}),$$
where
\begin{eqnarray*}
F_2 &=& -\int_0^{xy} \frac{\sqrt{2}e^{-t}(1+e^{-4t})}{\sqrt{2\pi}(1-e^{-4t})^{\frac{3}{2}}} (xy)^{-\frac{\a}{2}}(x-y) \exp\(-\frac{1+e^{-4t}}{1-e^{-4t}} \frac{(x-y)^2}{2}\)
\exp\(-\frac{(1-e^{-2t})^2}{1-e^{-4t}}xy\) \dtt , \cr
F_3 &=& -\int_0^{xy} \frac{1}{4 \sqrt{\pi} t^2} (xy)^{-\a\sl 2}(x-y) \exp\(- \frac{(x-y)^2}{4t}\)dt,\cr
F_4 &=& -\int_0^{\8} \frac{1}{4 \sqrt{\pi} t^2} (xy)^{-\a\sl 2}(x-y) \exp\(- \frac{(x-y)^2}{4t}\) dt = -\frac{(xy)^{-\a/2}}{\sqrt{\pi}(x-y)}.
\end{eqnarray*}
Similar analysis to that we have done in \eqref{est13} leads to
\begin{align}\label{est14}
|F_i-F_{i+1}|\leq C x^{-\a-1}, \ i=1,\dots,4.
\end{align}
Thanks to \eqref{est9}--\eqref{est14}, we have
\begin{equation}\label{all61}
\sup_{y<3} \int_{0}^{3} \chi_{\{|x-y|<y/2\}} \Big| R(x,y)-\frac{B}{x^{\a+1}-y^{\a+1}}- \frac{A-B}{x^{\a+1}+y^{\a+1}} \Big| d\mu(x) < \8.
\end{equation}
Observe that
\begin{equation}\label{all62}
\sup_{y<3} \int_{0}^{3}\chi_{\{|x-y|<y/2\}}\chi_{\{|x-y|>1/2\}} \(\Big|\frac{B}{x^{\a+1}-y^{\a+1}} \Big|+\Big|\frac{A-B}{x^{\a+1}+y^{\a+1}}\Big|\) d\mu(x) < \8.
\end{equation}
Note that if $|x-y|<1/2$ then $\phi\((x-y)/\rho(y)\)=1$ (see Figure 1). Therefore, it is not difficult to see that \eqref{all61}-\eqref{all62} imply
\begin{equation}\label{all6}
\sup_{y<3} \int_{0}^{3} \chi_{\{|x-y|<y/2\}}| g(x,y)| \dtt d\mu(x) < \8.
\end{equation}

Finally, the required estimate \eqref{integr} follows directly from \eqref{all1}, \eqref{all2}, \eqref{all3}, \eqref{all4}, \eqref{all5}, \eqref{all6}.
\end{proof}

{\bf Acknowledgments:} The author would like to thank Jacek Dziuba\'nski, Adam Nowak, and Krzysztof Stempak for their helpful comments and suggestions. The author is also greatly indebted to the referees for valuable remarks that improved the presentation of the paper.


\end{document}